\documentstyle[12pt]{article}

\newlength{\abstractwidth}
\renewcommand{\thefootnote}{\fnsymbol{footnote}}
\renewcommand{\thanks}[1]{\footnote{#1}} % Use this for footnotes
\newcommand{\starttext}{ \setcounter{footnote}{0}
\renewcommand{\thefootnote}{\arabic{footnote}}}

\newcommand{\be}{\begin{equation}}
\newcommand{\bea}{\begin{eqnarray}}
\newcommand{\eea}{\end{eqnarray}} 
\newcommand{\ee}{\end{equation}}

\def\ba{\begin{eqnarray}}
\def\ea{\end{eqnarray}}

\def\o{\omega}
\def\Re{{\rm Re}}

\def\det{{\rm det}}

\def\log{\,{\rm log}\,}

\def\o{\omega}

\def\o{\omega}

\def\p{\partial}

\def\ddb{{\partial\bar\partial}}

\def\[{{\bf [}}
\def\]{{\bf ]}}

\begin{document}
\starttext \baselineskip=15pt \setcounter{footnote}{0}
\newtheorem{theorem}{Theorem}
\newtheorem{lemma}{Lemma}
\newtheorem{definition}{Definition}
\newtheorem{proposition}{Proposition}
\newtheorem{corollary}{Corollary}
\newtheorem{remark}{Remark}

\begin{center}
{\Large \bf
Parabolic complex Monge-Amp\`ere equations on compact K\"ahler manifolds
\footnote{Contribution to the proceedings of the ICCM 2018 conference in Taipei, Taiwan. 
Work supported in part by the National Science Foundation Grant DMS-1809582.
Key words: complex Monge-Amp\`ere equations, geometric flows, concave equations, $C^k$ estimates, exponential convergence.
}
}
\end{center}

\centerline{Sebastien Picard and Xiangwen Zhang}

\bigskip

\begin{abstract}
We study the long-time existence and convergence of general parabolic complex Monge-Amp\`ere type equations whose second order operator is not necessarily convex or concave in the Hessian matrix of the unknown solution.
\end{abstract}

\medskip

\section{Introduction} 
Let $(X,\o)$ be a compact K\"ahler manifold. Let $f\in C^{\infty}(X, {\bf R})$ be a given function and $F: {\bf R}_+ \to {\bf R}$ a smooth strictly increasing function; that is $F'(\rho) >0$ for $\rho>0$.  In this paper, we consider the following parabolic complex Monge-Amp\`ere equation
\bea\label{flow}
\p_t u = F\left( e^{-f} \, \det ( \delta^i{}_j + \nabla^i \nabla_j u) \right),
\eea
with
\be
\o+ i\ddb u(x,t) >0.
\ee
Here we use the notation $\nabla^i \nabla_j u = g^{i \bar{k}} \nabla_{\bar{k}} \nabla_j u$, where $\omega = i g_{\bar{k} i} dz^i \wedge d \bar{z}^k$ and $g^{i\bar k}$ is the inverse matrix of $g_{\bar k i}$. One can check that $h^i{}_j = \delta^i{}_j + \nabla^i \nabla_j u$ defines an endomorphism of $T^{1,0} X$, and its determinant is a well-defined function on $X$. We study the long-time existence and convergence of the flow (\ref{flow}) with smooth initial data $u(x, 0) = u_0(x)$. We do not assume any concavity (or convexity) condition on the speed function $F$.                             

\smallskip

Complex Monge-Amp\`ere type equations have played an important role in the study of complex geometry because of the natural appearance of the determinant operator in the formula for the Chern-Ricci form ${\rm Ric} = - i\ddb \log \det \,\omega$ of a Hermitian metric $\omega$. The foundations of an existence and regularity theory for complex Monge-Amp\`ere equations on compact K\"ahler manifolds was pioneered by Yau \cite{Yau} (see also Aubin \cite{Aubin}), culminating in a complete solutions of the famous Calabi conjecture. Since then, this type of equation has been studied extensively and many new directions have emerged over the years. For example, a general theory of fully nonlinear elliptic concave equations was developed by Caffarelli, Kohn, Nirenberg and Spruck \cite{CNSI, CNSII, CNSIII}, a theory of generalized solutions for the complex Monge-Amp\`ere equation - or pluripotential theory - was laid out by Bedford-Taylor \cite{BT1, BT2}, and fundamental estimates in pluripotential theory were contributed by Kolodziej \cite{K1, K2, K3}. Some other works on this topic include \cite{BBGZ, BEGZ, DZZ, EGZ, FY1, FY2, GuanLi, Guan, STW, TW1, TW2, TW3, ZZ, Zhang}. We refer the interested reader to the survey paper on complex Monge-Amp\`ere equations by Phong-Song-Sturm \cite{PSS}.

%ellptic: metion Yau's work and motivate the remarkable development on general theory of concave fully nonlinear elliptic equations by CNS. See the nice survey paper on complex Monge-Ampere equaitons by Phong

There are also parabolic versions of the complex Monge-Amp\`ere equation which arise in complex geometry. In 1985, motivated by the Ricci flow initiated by Hamilton \cite{Hamilton}, H.-D. Cao \cite{Cao} studied the K\"ahler-Ricci flow on a compact K\"ahler manifold $(X, \o)$. In this setting, the Ricci flow reduces to the following parabolic complex Monge-Amp\`ere equation 
\bea\label{KRF}
\p_t u = \log \, \det (\delta^i{}_j + \nabla^i \nabla_j u) - f,
\eea
with $\o+ i\ddb u >0$. By adapting Yau's estimates to the parabolic setting, Cao obtained the long-time existence and convergence of this flow, and therefore provided a parabolic approach to Yau's solution of the Calabi conjecture. More recent parabolic proofs using the inverse Monge-Amp\`ere flow were obtained in \cite{CK, CHT,FLM}, which is of the form 
\bea\label{IMAF}
\p_t u = 1 - {e^f \over \det (\delta^i{}_j + \nabla^i \nabla_j u)}.
\eea
The inverse Monge-Amp\`ere flow also arises as the T-dual of the Anomaly flow \cite{FPic}. The Anomaly flows were introduced by the authors jointly with D.H. Phong \cite{PPZ2,PPZ5}. This family of geometric flows has been successful in investigating various facets of non-K\"ahler complex geometry and optimal metrics in theoretical physics \cite{FHP,FP,FPic,PPZ2,PPZ5,PPZ6,PPZ7,PPZ9,PPZ10}. In the special case of conformally K\"ahler initial data, the Anomaly flow reduces \cite{PPZ10} to a parabolic equation of the form
\bea \label{MAF}
 \p_t u = e^{-f} \det (\delta^i{}_j + \nabla^i \nabla_j u).
\eea
We will review the link between the Anomaly flow and this Monge-Amp\`ere equation in Section 2. 

\par In this paper, we study the long-time existence and convergence of the general parabolic Monge-Amp\`ere flow (\ref{flow}). It includes the K\"ahler-Ricci flow, the inverse Monge-Amp\`ere flow, and the conformally K\"ahler Anomaly flow as special cases, by simply taking $F(\rho)= \log \rho$, $F(\rho)= 1 - \rho^{-1}$, and $F(\rho)=\rho$. The theorem that we will prove is following:

\begin{theorem}
Let $(X, \o)$ be a compact K\"ahler manifold of complex dimension $n$. Let $f\in C^\infty(X, {\bf R})$ and $F\, : \, {\bf R}_+ \to {\bf R}$ be a strictly increasing smooth function. Then there exists a smooth solution $u$ to the parabolic complex Monge-Amp\`ere equation (\ref{flow}) for all time. Moreover, $\varphi = u - \int_X \, u\, \o^n$ converges in $C^\infty$ to a smooth function $\varphi_{\infty}$ satisfying
\be
(\omega + i \ddb \varphi_\infty)^n = c_0 \, e^f \omega^n, \ \ c_0 = {\int_X \omega^n \over \int_X e^f \omega^n}.
\ee
\end{theorem}

%parabolic: motivated by Hamilton, Cao gave a parabolic proof of Yau's theorem, the key is to study the following parabolic Monge-Ampere equation...... More recent parabolic proofs using inverse Monge- Amp`ere flow.....

%in this paper, we establish the long time exsitence and convergences of a general parabolic Monge-Ampere flow, which includes the K\"ahler-Ricci flow and inverse Monge-Ampere flow as special cases, by simply letting $f=0$ and taking $F(\rho)= \log \rho$ and $F(\rho)= - {1\over \rho}$, respectively. The theorem we will prove is as following:

%Theorem
In the above theorem, the function $F(\rho)$ is required to satisfy the strictly increasing condition which amounts to the parabolicity of the flow. In comparison with the well-studied parabolic complex Monge-Amp\`ere type flows such as (\ref{KRF}) and (\ref{IMAF}), an important new feature is that {\it we do not assume that the operator on the right-hand side is concave (or convex) in the Hessian matrix $i\ddb u$}. This is relevant in the analysis of the special case of the Anomaly flow (\ref{MAF}). Though the right-hand side of (\ref{MAF}), $\det(\delta^i{}_j + \nabla^i \nabla_j u)$, is perhaps the most natural and simplest complex Monge-Amp\`ere operator, unlike $\log \, \det$ or  $1/\det$, it is not concave in the Hessian. Therefore, it does not fall within the scope of standard PDE methods such as those developed in \cite{CNSIII, Guanbo, PT, Gabor}. In joint work with D.H .Phong \cite{PPZ10}, we overcame this difficulty in the conformally K\"ahler Anomaly flow case by introducing a new auxilary function for the $C^2$ estimate and developing new tools to obtain all the necessary estimates. As a byproduct, the Anomaly flow can be used to provide another proof of the classical theorem of Yau \cite{Yau} on the existence of Ricci-flat K\"ahler metrics. 

% In comparison with the well-studied parabolic complex Monge-Ampere type flows, the upshot of (?) is that NO concavity is assumed for the....Indeed, an particular case of (?) appeared naturally in the study of the Anomaly flow for conformal K\"ahler initial data. As we will see in section 2, it corresponds to the following simple equation
%u_t = det...

% This next paragraph is good, but the flows in this paper always converge so there is no stability condition to detect. May be best to leave it out and focus on the PDE aspects?

%As stressed in \cite{CHT} and also discussed in \cite{PPZ9}, even when the solution of an elliptic equation can already be found by a particular flow (e.g. a stability condition), other flows can be potentially useful too. This is because, in the absence of a stability condition, they may fail to converge in different ways, which would provide then different ways of detecting instability. 

 %Although different parabolic equations might lead to the same stationary points (that is their corresponding elliptic equations are equivalent), it is still interesting to study parabolic flows of different types. For the sake of geometry, as stressted in stressed in [14] and [PPZ], alternative approaches are interesting not just for alternative proofs 2 in themselves, but also for the study of singularities that they would develop when no stationary point exists. 

\par

From the PDE point of view, it is also of great importance to investigate different types of parabolic equations with possibly equivalent stationary points, as it will enrich our understanding of the structure of geometric PDEs. We note that the speed function can affect the dynamics and analysis of the flow. For a recent example of this in complex geometry, see the work of Fei-Guo-Phong \cite{FGP} on the LYZ flow \cite{LYZ}.
\par

To solve geometric problems, we are often forced to develop new tools to study fully nonlinear PDE with new features. For example, Guan \cite{Guanbo} and Sz\'ekelyhidi \cite{Gabor} introduced $\mathcal C$-subsolutions of fully nonlinear elliptic concave equations to develop a theory of solvability conditions for a large family of equations on compact manifolds, and Phong-T\^o \cite{PT} developed the parallel framework for parabolic concave equations. In \cite{PPZ6}, together with D.H. Phong, we developed some new techniques to obtain convergence of the Anomaly flow with Fu-Yau ansatz, and a similar idea was also used to solve the Fu-Yau Hessian equations \cite{PPZ8} in any dimension.
% (see also \cite{CHZ, FGV} for other applications of this idea). 
In the current paper, we generalize an idea from our earlier work with D.H. Phong \cite{PPZ10} on the Anomaly flow to study parabolic Monge-Amp\`ere flows without concavity.% We hope the new auxiliary function and idea might shed some lights on other geometric equations.

\par
Finally, we would like to mention that the flow (\ref{flow}) can be viewed as a complex version of the inverse Gauss curvature flow studied extensively in convex geometry, see for example \cite{CT, CT1, Tsai, Urbas}. In particular, Chow and Tsai studied the nonhomogeneous inverse Gauss curvature flows 
$\p_t u= F(\det\,  (u\, \delta_{ij} + \nabla_i \nabla_j u))$
on ${\bf S}^n$ under certain concavity assumption on the speed $F$. It might be worth pointing out that, due to the difference between complex and real cases, and also the absence of zeroth order term $u$ in the equation, we are able to prove our main theorem without any concavity assumption on the speed function. 
%
% the solvability of general parabolic concave equations under C-subsolution condition developed by Phong-To, which is a parabolic.... More recently, together with Phong, we also devolop some estimates with scales in the study of the Anomaly flow on Fu-Yau ansatz and also its elliptic extension to solve the Fu-Yau Hessian equations, 
%%

\par
The paper is organized as follows. In \S 2, we discuss several special cases of the flow (\ref{flow}), which naturally arise from the Anomaly flow in the study of the Hull-Strominger system. In \S 3, we derive all the evolution equations that will be subsequently needed. In \S 4, we establish all the necessary estimates for the long-time behavior of the flows and discuss convergence.

\medskip

\noindent
{\bf Acknowledgements}: The authors would like to thank their collaborator D.H. Phong for the wonderful collaborations and the inspiring discussions which led to this work. They are also indebted to him for his constant encouragement and support.

\

\section{Motivation from the Anomaly flow}
\setcounter{equation}{0}

As mentioned in the Introduction, some special cases of the flow (\ref{flow}) are closely related to the Anomaly flow with conformally K\"ahler initial data. In this section, we will briefly recall the reduction of the parabolic Monge-Amp\`ere equation (\ref{MAF}) from \cite{PPZ10}, which corresponds to the case $F(\rho)=\rho$ in (\ref{flow}). Moreover, we will also see that a slight modification of the Anomaly flow leads to the case $F(\rho) = \rho^a$.
%try to emphasize a point that: one of the hope of Anomaly flow is to develop new analytic techniques for PDE arising from geometric problems.

\medskip

\subsection{The Anomaly flow}
The Anomaly flow was introduced by the authors together with D.H. Phong in \cite{PPZ2} as a parabolic system incorporating the anomaly cancellation equation in the Hull-Strominger system \cite{Hull, Strominger}. More precisely, 
let $X$ be a compact complex Calabi-Yau $3$-fold, equipped with a nowhere vanishing holomorphic $(3, 0)$-form $\Omega$. Let $t\to \Phi(t)$ be a given path of closed $(2, 2)$-forms.
%, with $[\Phi(t)]=[\Phi_2(0)]$ for each $t$. 
Let $\omega_0$ be an initial metric which is conformally balanced \cite{LY}, meaning that
\be
d( \|\Omega\|_{\omega_0}\omega_0^2) =0,
\ee
where
$\| \Omega \|_\omega^2 \, {\omega^3 \over 3!} = i \Omega \wedge \bar{\Omega}$.
Then the Anomaly flow is the flow of $(2,2)$-forms defined by
\bea
\label{af1}
\partial_t(\|\Omega\|_\omega\omega^2)
=
i\partial\bar\partial\omega-{\alpha'\over 4}
({\rm Tr}\,(Rm\wedge Rm)-\Phi),
\eea
where $Rm$ is the curvature of the Chern unitary connection defined by the Hermitian metric $\omega$, viewed as a section of $\Lambda^{1,1}(X) \otimes End(T^{1,0}(X))$, and $\alpha'\in {\bf R}$ is a constant, called the slope parameter. We note that the stationary points of the flow are given by the anomaly cancellation equation in the Hull-Strominger system \cite{Hull,Strominger}, when $\Phi = {\rm Tr} \, F \wedge F$ and $F$ is the curvature of a connection on a gauge bundle $E \rightarrow X$.

\medskip

An important feature of the Anomaly flow is that {\it it preserves the conformally balanced condition} of the initial data. Indeed, by Chern-Weil theory, its right-hand side is a closed $(2,2)$-form and hence for all $t$,
$$
d(\|\Omega\|_\omega\omega^2)=0
$$
if $d(\|\Omega\|_\omega\omega^2)=0$ at time $t=0$.

In a sequence of papers \cite{PPZ2, PPZ5, PPZ6, PPZ7, PPZ10}, joint with D.H. Phong, the authors studied the behavior of the Anomaly flow from different aspects. For example, in \cite{PPZ5}, we rewrote the flow of $(2, 2)$-forms as a flow of the Hermitian metric $\o$, and this local expression allowed us to partially carry out the Shi-type curvature estimates for the derivatives of the curvature and torsion tensors, which is an important step toward the analysis of the long-time behavior of the flow. In \cite{PPZ7}, we considered the Anomaly flow on toric fibrations over a $K3$ surface, which is the geometric setting of the first non-K\"ahler solution of the Hull-Strominger system found by Fu-Yau \cite{FY1, FY2} (for recent work on the Fu-Yau compactifications, see also e.g. \cite{FGV,GF,PPZ1,PPZ4} and references therein). By developing new analytic tools, we were able to prove the long-time existence and convergence of the flow in this setting, and therefore provide an alternative proof of Fu-Yau's result which works for all $\alpha'\in {\bf R}$. 

However, due to several substantial new difficulties such as the non-vanishing torsion and higher order curvature terms, a full understanding of the Anomaly flow is still out of reach at this moment. Therefore, it is our hope to develop more robust analytic tools to study this system of equations as we move toward the goal of obtaining a criteria for long-time existence and convergence. 

%some natural geometric PDEs but without the properties that required to apply the general fully nonlinear PDE theory. 
\medskip

\subsection{The case $F(\rho)=\rho$}

To obtain a better understanding for (\ref{af1}), we studied a model case of the Anomaly flow in \cite{PPZ10} by taking $\alpha'=0$. On an $n$-dimensional compact complex manifold $X$ equipped with a nonvanishing holomorphic $(n, 0)$- form $\Omega$, we consider the flow
\bea\label{af2}
\p_t \left(\|\Omega\|_{\omega} \, \omega^{n-1} \right) = i\ddb \omega^{n-2}
\eea
with $\omega(0) = \omega_0$ satisfying the conformally balanced condition 
\bea\label{cbc}
d\left(\|\Omega\|_{\omega} \, \omega^{n-1}\right) =0.
\eea
 It is easy to see that this reduces to the zero slope case of the flow (\ref{af1}) when we let the dimension be equal to $n=3$.

In fact, other than serving as a model case for the general Anomaly flow (\ref{af1}), the flow (\ref{af2}) is also interesting on its own as a non-K\"ahler flow preserving the conformally balanced condition. If the flow does converge to a limiting metric $\omega_\infty$, then $\omega_\infty$ must be both conformally balanced $d\left(\|\Omega\|_{\omega_\infty} \, \omega_\infty^{n-1}\right) =0$ and astheno-K\"ahler $i\ddb \o_{\infty}^{n-2}=0$. Manipulating these identities and applying the maximum principle (see \cite{FT, MT, PPZ10}) shows that $\omega_\infty$ must be a K\"ahler Ricci-flat metric. Therefore, the flow (\ref{af2}) can in principle be used to determine whether a conformally balanced manifold is actually K\"ahler.

\medskip

In \cite{PPZ10}, we establish a partial converse of the above observation, which states that if we start from a conformal K\"ahler initial metric satisfying the conformal balanced condition, then the flow (\ref{af2}) will preserve the conformal K\"ahler ansatz, exist for all time and converge to a K\"ahler Ricci-flat metric. More precisely, suppose $X$ is a compact K\"ahler manifold and let $\omega_0$ be a Hermitian metric satisfying
\bea\label{omega0}
\|\Omega\|_{\omega_0} \, \omega_0^{n-1} = \chi^{n-1}
\eea
for some K\"ahler metric $\chi$. Equivalently, this means that $\omega_0$ is conformally K\"ahler, as $\o_0=\|\Omega\|^{-2/n-2}_{\chi} \, \chi$. Now, we take $t\to u(t)$ to be the solution of the solution of the Monge-Amp\`ere flow 
\bea\label{MAf}
\p_t u= e^{-f} \, {\det(\chi + i\ddb u) \over \det \, \chi}
\eea
with $f\in C^\infty(X, {\bf R})$ given by 
\bea
e^{-f} = {1\over n-1} \, \|\Omega\|_{\chi}^{-2}.
\eea
To connect the flow (\ref{af2}) with the above Monge-Amp\`ere flow, we consider Hermitian metrics $t\to \omega(t)$ satisfying
\bea\label{ansatz}
\|\Omega\|_{\omega(t)}\, \omega^{n-1}(t) = (\chi+i \ddb u(t))^{n-1}
\eea
with $u(t)$ being the solution of the flow (\ref{MAf}). It is easy to see from (\ref{omega0}) that, when $t=0$, $\omega_0$ satisfies this ansatz with $u(\cdot\,, 0)=0$. In fact, the calculation in \cite{PPZ10} shows that $\omega(t)$ given by (\ref{ansatz}) also satisfies the evolution equation (\ref{af2}) for any $t>0$. Now, from the strict parabolicity and uniqueness of the solution of flow (\ref{af2}), we conclude that the long time behavior of (\ref{af2}) with initial data (\ref{omega0}) reduces to the study of the Monge-Amp\`ere flow (\ref{MAf}) with $u(\cdot\, , 0)=0$. This is exactly the case when $F(\rho)=\rho$ in our main equation (\ref{flow}).

\medskip

\subsection{The case $F(\rho)= \rho^a$}

In this subsection, we will see that the flow (\ref{flow}) with $F(\rho) = \rho^a$ also has a natural connection to a modified version of the Anomaly flow with zero slope parameter.

\smallskip

For this, we simply generalize the definition for the conformally balanced condition as
\bea\label{gcbc}
d\left(\|\Omega\|^\beta_{\omega} \, \omega^{n-1}\right) =0,
\eea
for $\beta \in {\bf R}$, and then we modify the flow (\ref{af2}) as
\bea\label{newflow}
\p_t \left(\|\Omega\|^\beta_{\omega} \, \omega^{n-1} \right) = (-1)^\sigma\, i\ddb \omega^{n-2}
\eea
with $\sigma\in \{0,1\}$ and
$\omega(0) = \omega_0$ satisfying (\ref{gcbc}). It is obvious that the modified conformally balanced condition is still preserved by the modified flow (\ref{newflow}).  The case when $\beta = n-1$ and $\lambda =1$, given by
\be \label{dual-af}
\p_t  ( \| \Omega \|_\omega \omega)^{n-1} = - i \ddb \omega^{n-2},
\ee
is of particular interest. It was first derived by Fei-Picard \cite{FP} in the context of mirror symmetry, and it arises as the evolution equation induced by the Anomaly flow on the dual torus fibration of semi-flat Calabi-Yau $n$-folds.

%We note that the choice of $\sigma$ depends on $\beta$ which is purely related to the parabolicity of the flow (\ref{newflow}). One can see this by computing the symbol of the operator. We would rather skip the detailed computation here since the flow (\ref{newflow}) is not the main concern of this paper. But, there are two important cases providing evidence for this: one is the Anomaly flow (\ref{af2}) discussed in the previous section and studied in detail in \cite{PPZ10}, for which $\beta=1$ and $\sigma=0$; the other interesting case is when $\beta=n-1$ and $\sigma=1$, with the flow given by
%\be \label{dual-af}
%\p_t  ( \| \Omega \|_\omega \omega)^{n-1} = - i \ddb \omega^{n-2}.
%\ee
%This flow was derived by Fei-Picard \cite{FP} in the context of mirror symmetry, and it arises as the evolution equation induced by the Anomaly flow (\ref{af2}) on the dual torus fibration of semi-flat Calabi-Yau $n$-folds. 

\smallskip

Now, following the same discussion as in previous section, we study the modified flow (\ref{newflow}) with conformally K\"ahler initial data. Suppose the initial metric satisfies 
\bea
\|\Omega\|^\beta_{\omega_0} \, \omega_0^{n-1}= \chi^{n-1}
\eea 
for some K\"ahler metric $\chi$. Then, by the same calculation, one can show that the Hermitian metrics $t\to \omega(t)$ defined by 
\bea
\|\Omega\|^\beta_{\omega(t)} \, \omega(t)^{n-1} = (\chi + i\ddb u(t))^{n-1}, \ \ {\rm with} \ \chi+i \ddb u>0
\eea
will satisfy the flow (\ref{newflow}) if $t\to u(t)$ is a solution of the following parabolic Monge-Amp\`ere flow
\bea\label{para-equ}
\p_t u = (-1)^\sigma {1\over n-1}\left(e^{-f} \cdot {\det(\chi+i\ddb u)\over \det \chi} \right)^a,
\eea
with $e^{-f} = \|\Omega\|_{\chi}^{-2}$ and $a={(n-2)\beta\over 2n-2-n\beta}$. This is the case when $F(\rho)= \pm \rho^a$ in the flow (\ref{flow}).

\smallskip
For the dual Anomaly flow (\ref{dual-af}), we have $\beta = n-1$ and $\sigma=1$, and we obtain
\be
\p_t u = -{1\over n-1}\left(e^{-f} \cdot {\det(\chi+i\ddb u)\over \det \chi} \right)^{-1}.
\ee
This coincides with the inverse Monge-Amp\`ere flow (\ref{IMAF}) which was previously studied in K\"ahler geometry by Cao-Keller \cite{CK}. It is the vanishing first Chern class case of the $MA^{-1}$ flow introduced by Collins-Hisamoto-Takahashi \cite{CHT}, which is the gradient flow of the Ding energy functional, and is used to study the connection between K\"ahler-Einstein metrics and $K$-stability.

\

\section{Evolution equations}
\setcounter{equation}{0}

Let $(X,\chi)$ be a compact K\"ahler manifold. In this section, we compute the evolution of various quantities along the flow of $u(x,t)$ given by
\bea\label{FefMA}
\p_t u = F\left( e^{-f} \, \det \, h \right), \ \ 
\eea
with 
\be
h^i{}_j = \delta^i{}_j + \nabla^i \nabla_j u > 0,
\ee
and $f \in C^\infty(X,{\bf R})$ and $F: {\bf R}_+ \rightarrow {\bf R}$ are given functions. We require that $F$ satisfies the ellipticity condition $F'>0$.
\smallskip
\par Our conventions, fixed from this point on in the paper, are as follows. We fix a background metric $\chi = i \chi_{\bar{k} j} dz^j \wedge d \bar{z}^k$. We denote the inverse of $\chi_{\bar{k} j}$ by $\chi^{j \bar{k}}$, so that $\chi^{i \bar{k}} \chi_{\bar{k} j} = \delta^i{}_j$. All norms are with respect to $\chi$. The covariant derivative $\nabla$ is the Chern connection with respect to $\chi$. This is defined by acting on sections $W \in (T^{1,0}X)^*$, locally written as $W= W_i dz^i$, by
\be
\nabla_k W_i = \partial_k W_i - \Gamma^r_{ki} W_r, \ \ \nabla_{\bar{k}} W_i = \partial_{\bar{k}} W_i,
\ee
where $\Gamma^r_{ki} = \chi^{r \bar{p}} \partial_k \chi_{\bar{p} i}$. 
\smallskip
\par Curvature tensors $R_{\bar{k} j}{}^r{}_i = - \p_{\bar{k}} \Gamma^r_{ji}$ are with respect to the background $(\chi,\nabla)$. The Laplacian of $\chi$ on functions is denoted
\be
\Delta = \chi^{p \bar{q}} \nabla_p \nabla_{\bar{q}}.
\ee
We write
\be
g_{\bar{k} j} = \chi_{\bar{k} j} + u_{\bar{k} j},
\ee
for the evolving metric, and $g^{j \bar{k}}$ for the inverse $g^{-1}$. Then the endomorphism $h$ can also be written as
\be
h^i{}_j = \chi^{i \bar{k}} g_{\bar{k} j} = \chi^{i \bar{k}} (\chi_{\bar{k} j} + \nabla_{\bar{k}} \nabla_j u).
\ee
We will use the linearized operator
\be \label{L-defn}
L = F' e^{-f} \det \, h \, g^{j\bar k}\, \nabla_j \nabla_{\bar k}.
\ee

\subsection{Evolution of $u$}

We compute
\be
(\p_t - L) u = F - F' e^{-f} \, \det \, h \, g^{j\bar k} \nabla_j\nabla_{\bar k} u.
\ee
Since $u_{\bar{k} j} = g_{\bar{k} j} - \chi_{\bar{k} j}$,
\be \label{evol-u}
(\p_t - L) u = F - n F' e^{-f} \det \, h + F' e^{-f} \det \, h \, g^{j\bar k} \chi_{\bar{k} j}.
\ee
The normalized function $\varphi = u - V^{-1} \int_X u$, where $V= \int_X \chi^n$, evolves by
\be \label{evol-phi}
(\p_t - L) \varphi = F - n F' e^{-f} \det \, h + F' e^{-f} \det \, h \, g^{j\bar k} \chi_{\bar{k} j} - {1 \over V} \int_X F \, \chi^n.
\ee

\subsection{Evolution of $F$}
The equation $\p_t u = F$, differentiated in time is
\be
{d \over dt} \p_t u = L \p_t u,
\ee
thus
\be
(\p_t - L) F = 0,
\ee
and
\be
(\p_t - L) F^2 = - 2 F' e^{-f} \, \det \, h\, g^{j\bar k} \nabla_j F \nabla_{\bar{k}} F.
\ee
Then
\be
(\p_t - L) F^2 = - 2 (F')^3 e^{-f} \, \det \, h\, g^{j\bar k} \nabla_j (e^{-f} \det \, h) \nabla_{\bar{k}} (e^{-f}\det \, h ).
\ee
Distributing the derivatives, we obtain
\bea \label{evol-F2}
(\p_t - L) F^2 &=&- 2 (F')^3 e^{-3f} \, \det \, h\, g^{j\bar k} \nabla_j \det \, h  \nabla_{\bar{k}} \det \, h \nonumber\\
&&- 2 (F')^3 e^{-3f} \, (\det \, h)^3 g^{j\bar k} \nabla_j f \nabla_{\bar{k}} f \nonumber\\
&&+ 4 (F')^3 e^{-3f} \, (\det \, h)^2 \, \Re \{ g^{j\bar k} \nabla_j f \nabla_{\bar{k}} \det \, h \}.
\eea

\subsection{Evolution of $h$}
Differentiating the equation $\p_t u = F$ gives
\be
{d \over dt} \nabla_j u = F' \nabla_j (e^{-f} \det \, h).
\ee
Differentiating again gives
\be
{d \over dt} \Delta u = F' \Delta (e^{-f} \det \, h) +  F'' |\nabla (e^{-f} \det \, h)|^2.
\ee
Distributing derivatives, we obtain
\bea \label{evol-Delta}
    {d \over dt} \Delta u &=& F' e^{-f} \Delta \det \, h - 2 F' e^{-f} \Re \langle \nabla f, \nabla \det \, h \rangle \nonumber\\
    &&+  F'' |\nabla (e^{-f} \det \, h)|^2 + F' \det \, h \Delta e^{-f}.
\eea
Using $\delta \log \det \, h = {\rm Tr} \, h^{-1} \delta h$, we differentiate $\det \, h$ once and get
\be
\nabla_{\bar{q}} \det \, h = (\det \, h) {\rm Tr} \, h^{-1} \nabla_{\bar{q}} h.
\ee
Using $\delta h^{-1} = - h^{-1} \,\delta h\, h^{-1}$, we differentiate again and get
\bea
\nabla_p \nabla_{\bar{q}} \det \, h &=& (\det \, h) {\rm Tr} \, h^{-1} \nabla_p \nabla_{\bar{q}} h - (\det \, h) {\rm Tr} \, h^{-1} \nabla_p h  h^{-1} \nabla_{\bar{q}} h \nonumber\\
&&+ {1 \over \det \, h} \nabla_p \det \, h \nabla_{\bar{q}} \det \, h .
\eea
Since $h^i{}_j = \chi^{i \bar{k}} g_{\bar{k} j}$, $(h^{-1}){}^i{}_j= g^{i \bar{k}} \chi_{\bar{k} j}$, and $\nabla \chi_{\bar k j} = 0$, we obtain
\bea \label{Delta-deth}
\Delta \det \, h &=& (\det \, h) \chi^{p \bar{q}} g^{j \bar{k}} \nabla_p \nabla_{\bar{q}} \nabla_j \nabla_{\bar{k}} u - (\det \, h) \chi^{p \bar{q}} g^{j \bar{r}} g^{s \bar{k}} \nabla_p g_{\bar{r} s}  \nabla_{\bar{q}} g_{\bar{k} j} \nonumber\\
&&+ {1 \over \det \, h} |\nabla \det \, h|^2.
\eea
We use the convention, for any section $W = W_i dz^i$ of $(T^{1,0} X)^*$,
\be
(\nabla_j \nabla_{\bar{k}}-\nabla_{\bar{k}} \nabla_j) W_i = -R_{\bar{k} j}{}^r{}_i W_r.
\ee
Exchanging covariant derivatives gives
\be
\nabla_p \nabla_{\bar{q}} \nabla_j \nabla_{\bar{k}} u = \nabla_p \nabla_{\bar{k}} \nabla_j \nabla_{\bar{q}} u =  \nabla_{\bar{k}} \nabla_p  \nabla_j \nabla_{\bar{q}} u - R_{\bar{k} p}{}^r{}_j u_{\bar{q} r} + R_{\bar{k} p \bar{q}}{}^{\bar{r}} u_{\bar{r} j}.
\ee
Thus
\be
\chi^{p \bar{q}} \nabla_p \nabla_{\bar{q}} \nabla_j \nabla_{\bar{k}} u=  \nabla_{\bar{k}}  \nabla_j \Delta u - \chi^{p \bar{q}} R_{\bar{k} p}{}^r{}_j u_{\bar{q} r} + R_{\bar{k}}{}^{\bar{r}} u_{\bar{r} j}.
\ee
Writing $u_{\bar{q}r} = g_{\bar{q}r} - \chi_{\bar{q} r}$ and contracting with $g^{j \bar{k}}$, this simplifies to
\be
g^{j \bar{k}} \chi^{p \bar{q}} \nabla_p \nabla_{\bar{q}} \nabla_j \nabla_{\bar{k}} u =  g^{j \bar{k}} \nabla_{\bar{k}}  \nabla_j \Delta u - g^{j \bar{k}} R_{\bar{k} p}{}^r{}_j h^p{}_r + R.
\ee
Substituting this expression into (\ref{Delta-deth}) gives
\bea
\Delta \det \, h &=& (\det \, h)  g^{j \bar{k}} \nabla_{\bar{k}}  \nabla_j \Delta u - (\det \, h) \chi^{p \bar{q}} g^{j \bar{r}} g^{s \bar{k}} \nabla_p g_{\bar{r} s}  \nabla_{\bar{q}} g_{\bar{k} j} \nonumber\\
&&- (\det \, h) g^{j \bar{k}} R_{\bar{k} p}{}^r{}_j h^p{}_r +  (\det \, h) R + {1 \over \det \, h} |\nabla \det \, h|^2.
\eea
Substituting this into (\ref{evol-Delta}), we obtain
\bea \label{evol-Delta2}
    {d \over dt} \Delta u &=& F' e^{-f}(\det \, h) \bigg\{   g^{j \bar{k}} \nabla_{\bar{k}}  \nabla_j \Delta u - \chi^{p \bar{q}} g^{j \bar{r}} g^{s \bar{k}} \nabla_p g_{\bar{r} s}  \nabla_{\bar{q}} g_{\bar{k} j} \nonumber\\
&&-  g^{j \bar{k}} R_{\bar{k} p}{}^r{}_j h^p{}_r +   R + {1 \over (\det \, h)^2} |\nabla \det \, h|^2 \nonumber\\
    && - {2 \over \det \, h}  \Re \langle \nabla f, \nabla \det \, h \rangle +  {F'' e^f \over F' (\det \, h)} |\nabla (e^{-f} \det \, h)|^2 +  e^f \Delta e^{-f}\bigg\} .
    \eea
Recall the definition of $L$ (\ref{L-defn}). A direct calculation gives
    \be
(\p_t - L) \log {\rm Tr} \, h = {1 \over {\rm Tr} \, h} (\p_t - L) {\rm Tr} \, h + {F' e^{-f}(\det \, h) \over ({\rm Tr} \, h)^2} g^{j \bar{k}} \nabla_{\bar{k}}{\rm Tr} \, h  \nabla_j {\rm Tr} \, h.
    \ee
    Since ${\rm Tr} \, h = n + \Delta u$, substituting (\ref{evol-Delta2}) yields
    \bea \label{evol-Trh}
    (\p_t - L) \log {\rm Tr} \, h &=& {F' e^{-f}(\det \, h) \over {\rm Tr} \, h} \bigg\{  {1 \over {\rm Tr} \, h} g^{j \bar{k}} \nabla_{\bar{k}}{\rm Tr} \, h  \nabla_j {\rm Tr} \, h - \chi^{p \bar{q}} g^{j \bar{r}} g^{s \bar{k}} \nabla_p g_{\bar{r} s}  \nabla_{\bar{q}} g_{\bar{k} j} \nonumber\\
&&-  g^{j \bar{k}} R_{\bar{k} p}{}^r{}_j h^p{}_r +   R  + {1 \over (\det \, h)^2} |\nabla \det \, h|^2 - {2 \over \det \, h}  \Re \langle \nabla f, \nabla \det \, h \rangle \nonumber\\
&&+  {F'' e^{-f} |\nabla f|^2 \det \, h \over F'} +  {F'' e^{-f} \over F' (\det \, h)} |\nabla \det \, h|^2 \nonumber\\
&&-  2 {F'' e^{-f} \over F'} \langle \nabla f, \nabla \det \, h \rangle +  e^f\Delta e^{-f} \bigg\} .
    \eea
    We note that in the case of the K\"ahler-Ricci flow, then $F(\rho)=\log \rho$ and the function $f$ is the Ricci potential of the background metric $\chi$ (i.e.  $\partial_j \partial_{\bar{k}} f = R_{\bar{k} j}$), and the equation simplifies to
 \bea
 (\p_t - g^{j \bar{k}} \nabla_j \nabla_{\bar{k}}) \log {\rm Tr} \, h &=& {1 \over {\rm Tr} \, h} \bigg\{ -  g^{j \bar{k}} R_{\bar{k} p}{}^r{}_j h^p{}_r + {1 \over {\rm Tr} \, h} g^{j \bar{k}} \nabla_{\bar{k}}{\rm Tr} \, h  \nabla_j {\rm Tr} \, h  \nonumber\\
 &&- \chi^{p \bar{q}} g^{j \bar{r}} g^{s \bar{k}} \nabla_p g_{\bar{r} s}  \nabla_{\bar{q}} g_{\bar{k} j} \bigg\},
 \eea
which agrees with e.g. Equation (2.18) in Song-Weinkove \cite{SW}.

\

\section{Estimates}
\setcounter{equation}{0}
Using the calculations done in the previous section, we now work to obtain time-independent estimates along the flow. We will use the usual convention where constants $C$ may change line by line along the course of a proof, but only depend on the quantities stated in the theorem.

\subsection{Estimate of the determinant}

We consider the evolution equation for function $H=e^{-f} \, \det \, h$:
\bea
\p_t H&=& e^{-f} \, \det \, h \, g^{j\bar k} \nabla_j \nabla_{\bar k} \, \p_t u
\nonumber\\
&=&
e^{-f} \, \det \, h\, g^{j\bar k} \nabla_j\nabla_{\bar k} F(H)
\nonumber\\
&=&
e^{-f} \, \det \, h \, g^{j\bar k} \left(F' \, \nabla_j \nabla_{\bar k} H + F'' \, H_j\, H_{\bar k} \right).
\eea
By the maximum principle, we obtain
\bea
\p_t H_{\max} \leq 0, \ \ \ \p_t H_{\min} \geq 0.
\eea
It follows that
\bea
H_{\max}(t) \leq H_{\max }(0),\ \ \ H_{\min}(t) \geq H_{\min}(0),
\eea
and we have the following estimate:

\begin{lemma} \label{lemma-deth}
Let $u$ solve $\p_t u = F(e^{-f} \det \, h)$ with $\chi+ i \ddb u >0$ on $X \times [0,T]$. There exists a constant $C>0$ depending on $(X,\chi)$, $f$, and the initial data $u_0$, such that
\be
C^{-1} \leq \det \, h \leq C.
\ee
Consequently, we also have
\be
|\p_t u | \leq C,
\ee
where $C$ depends on $(X,\chi)$, $f$, $F$, and the initial data $u_0$.
\end{lemma}

\subsection{Uniform estimate}
Consider the function
\be
\varphi(x,t) = u(x,t) - {1 \over V} \int_X u \, \chi^n,
\ee
where $V = \int_X \chi^n$. In this section, we will prove an $L^\infty$ estimate for $\varphi$.

\begin{lemma} \label{lemma-c0}
Let $u$ solve $\p_t u =F(e^{-f} \det \, h)$ with $\chi+i \ddb u>0$ on $X \times [0,T]$. Let
\be
\varphi = u - {1 \over V} \int_X u \, \chi^n.
\ee
Then
\be
\| \varphi \|_{L^\infty(X \times [0,T])} \leq C,
\ee
where $C$ depends on the initial data, $(X,\chi)$, $f$, $F$.
\end{lemma}

{\it Proof:} Let
\be
w(x,t) = u(x,t) - (\sup_X u)(t).
\ee
We will show that at each $t \in [0,T]$, there holds
\be \label{w-lower}
(\inf_X w)(t) \geq -C.
\ee
This would imply the oscillation bound
\be \label{u-osc}
0 \leq (\sup_X u - \inf_X u)(t) \leq C,
\ee
at each time $t \in [0,T]$. By definition of $\varphi$, we have
\be
\int_X \varphi \, \chi^n = 0,
\ee
along the flow, and furthermore (\ref{u-osc}) implies $\sup_X \varphi - \inf_X \varphi \leq C$. This gives full control of $\| \varphi \|_{L^\infty}$, and Lemma \ref{lemma-c0} follows.
\bigskip
\par Thus we need to prove (\ref{w-lower}). We showed in the previous section that
\be \label{est-detw}
C^{-1} \leq {(\chi+ i \ddb w)^n \over \chi^n} \leq C
\ee
along the flow. The lemma is then a consequence of Yau's $C^0$ estimate \cite{Yau}. For completeness of exposition, we provide the full details following the argument due to Blocki \cite{Blocki} (see also Sz\'ekelyhidi \cite{Gabor} and Phong-T\^o \cite{PT}).
\bigskip
\par Before starting the estimate, we recall the definition of the contact set. For a smooth function $v: \overline{B}_1(0) \rightarrow {\bf R}$, we define
\be
\Gamma^+(v) = \{ x \in B_1(0) : v(x) \geq 0, \ \ v(y) \leq v(x) + Dv(x) \cdot(y-x), \ \ {\rm for} \ {\rm all} \ y \in B_1(0) \}.
\ee
For $\epsilon>0$, we also define
\be
\Gamma^+_\epsilon(v) = \Gamma^+(v) \cap \{ v \geq \sup_{B_1(0)} v - \epsilon \}.
\ee
We will use the following version of the Aleksandrov-Bakelman-Pucci (ABP) estimate of elliptic PDE theory.
\begin{proposition}
 Suppose $v: \overline{B}_1(0) \rightarrow {\bf R}$ is a smooth function with $\sup_{B_1} v = M > 0$. Let $\epsilon>0$. Suppose
\be
v|_{\partial B_1} \leq M - \epsilon.
\ee
Then
\be
\epsilon^n \leq C(n) \int_{\Gamma^+_\epsilon(v)} |\det D^2 v|.
\ee
\end{proposition}

{\it Proof:} Let $\psi = v - M + \epsilon$. Then $\sup_{B_1} \psi = \epsilon$ and $\sup_{\partial B_1} \psi \leq 0$. Applying the standard ABP estimate (e.g. Lemma 5.7 in \cite{HanLin}) to $\psi$ gives the result. Q.E.D.

\bigskip
\par Next, we note Blocki's adaptation of the ABP estimate to complex PDE.

\begin{proposition} \label{blocki-abp} \cite{Blocki} Suppose $v: \overline{B}_1(0) \subset {\bf C}^n \rightarrow {\bf R}$ is a smooth function with $\sup_{B_1} v = M > 0$. Let $\epsilon>0$. Suppose
\be
v|_{\partial B_1} \leq M - \epsilon.
\ee
Then
\be
\epsilon^n \leq C(n) \int_{\Gamma^+_\epsilon(v)} |\det v_{\bar{k} j}|^2.
\ee
\end{proposition}

{\it Proof:} On the set $\Gamma^+_\epsilon(v)$, the function $v$ is concave and its real Hessian satisfies
\be
D^2 v \leq 0.
\ee
Blocki \cite{Blocki} proved that for an arbitrary function $v: {\bf C}^n \rightarrow {\bf R}$, at points where $D^2 v \leq 0$, the pointwise inequality
\be
\det (-D^2 v) \leq 2^{2n} (\det v_{\bar{k} j})^2
\ee
holds. The proposition follows from the previous one. Q.E.D.
\bigskip
\par We now return to the estimate along the flow and the proof of Lemma \ref{lemma-c0}. We will follow the presentation of \cite{Gabor}. Fix $t_0 \in [0,T]$. We localize to a point $p \in X$ where $\inf_{X \times \{ t_0 \}} w$ is attained. Choose a coordinate chart $B_1(0)$ such that this point $p$ corresponds to the origin. In other words,
\be
w(0,t_0) = -M,
\ee
where $-M = \inf_{X \times \{ t_0 \}} w$ and $M \geq 0$. We may assume $M>0$, otherwise there is nothing to do as (\ref{w-lower}) holds and the proof of Lemma \ref{lemma-c0} is complete. Let $\varepsilon>0$. Following Sz\'ekelyhidi \cite{Gabor}, define the local function $v: B_1(0) \rightarrow {\bf R}$ by
\be
v(x) = -w(x,t_0) - \varepsilon |z|^2.
\ee
This function satisfies
\be
\sup_{B_1} v = M>0, \ \ \sup_{\partial B_1} v \leq M - \varepsilon.
\ee
Applying Blocki's version of the ABP estimate (Proposition \ref{blocki-abp}), we have
\be
\varepsilon \leq C \bigg( \int_{\Gamma_\varepsilon^+(v)} |\det v_{\bar{k} j}|^2 \bigg)^{1/n}.
\ee
On $\Gamma^+_\varepsilon(v)$, we have $D^2 v \leq 0$ and hence $D^2 w|_{t=t_0} \geq - \varepsilon I$. Choose $\varepsilon>0$ depending on $(X,\chi)$ small enough such that
\be
- \varepsilon \delta_{kj} \geq - {1 \over 2} \chi_{\bar{k} j}
\ee
on $B_1(0)$. Let $\lambda_1 \geq \cdots \geq \lambda_n > 0$ denote the eigenvalues of the endomorphism
\be
h^i{}_j |_{t=t_0} = \delta^i{}_j + \chi^{i \bar{k}}(w|_{t=t_0})_{\bar{k} j}.
\ee
Then on $\Gamma_\varepsilon^+(v)$, we have
\be
\lambda_n \geq {1 \over 2}.
\ee
By (\ref{est-detw}), we know
\be
\lambda_1 \lambda_2 \cdots \lambda_n \leq C.
\ee
Therefore
\be
\lambda_1 \leq 2^{n-1} \lambda_1 \lambda_2 \cdots \lambda_n \leq C,
\ee
and hence $|\det v_{\bar{k} j}| \leq C$. It follows that
\be
\varepsilon^n \leq C |\Gamma_\varepsilon^+(v)| .
\ee
However, by the definition of $\Gamma_\varepsilon^+(v)$, we have for any $p>0$,
\be
|M-\varepsilon|^p \leq {1 \over |\Gamma_\varepsilon^+(v)|}\int_{\Gamma_\varepsilon^+(v)} |v|^p,
\ee
and hence
\be \label{M-Lp-v}
|M- \varepsilon| \leq C \bigg( \int_{B_1} |v|^p \bigg)^{1/p}.
\ee
To control $-\inf_{X \times \{ t_0 \}} w = M$, it remains to control the $L^p$ norm of $v$ for some $p>0$. We have by definition
\be
|v(z)|^p = \left|-u(z,t_0) + \sup_{X \times \{t_0\}} u  - \varepsilon |z|^2 \right|^p,
\ee
and hence
\be \label{v-2-w}
|v|^p \leq \left| \sup_{X \times \{t_0\}}u - u(\cdot,t_0) \right|^p + C.
\ee
Since
\be
n+ \Delta u = \chi^{j \bar{k}} (\chi_{\bar{k} j} + u_{\bar{k} j}) > 0,
\ee
the function
\be
\tilde{w} = \sup_{X \times \{t_0\}}u - u(\cdot,t_0)
\ee
satisfies $\Delta \tilde{w} \leq n$. We also know that $\tilde{w} \geq 0$ and $\tilde{w}$ attains $0$ at a point $z \in X$. We can therefore apply the weak Harnack inequality (e.g. Theorem 4.15 in \cite{HanLin}) to obtain the existence of $p>0$ such that
\be \label{covering1}
\bigg( \int_{B_{1/2}(z)} |\tilde{w}|^p \bigg)^{1/p} \leq C(1+ \inf_{B_{1/2}(z)} \tilde{w}) \leq C,
\ee
where $B_1(z)$ is a local coordinate ball centered at $z$. Let $B_1(z_2)$ be another coordinate ball such that $B_{1/2}(z) \cap B_{1/2}(z_2) \neq \emptyset$. Then by the weak Harnack inequality,
\be
\bigg( \int_{B_{1/2}(z_2)} |\tilde{w}|^p \bigg)^{1/p} \leq C(1+ \inf_{B_{1/2}(z_2)} \tilde{w}) \leq C(1+ \inf_{B_{1/2}(z_2) \cap B_{1/2}(z)} \tilde{w}).
\ee
Thus
\be
\bigg( \int_{B_{1/2}(z_2)} |\tilde{w}|^p \bigg)^{1/p} \leq C \left( 1+ \bigg( \int_{B_{1/2}(z_2) \cap B_{1/2}(z)} |\tilde{w}|^p \bigg)^{1/p} \right) \leq C.
\ee
Covering $X$ by coordinate balls $B_{1/2}(z_i)$, we obtain
\be \label{covering2}
\int_X |\tilde{w}|^p \chi^n \leq C.
\ee
It follows from (\ref{v-2-w}) that
\be
\int_{B_1} |v|^p \leq C,
\ee
and from (\ref{M-Lp-v}) that $\inf_{X \times \{t_0 \}} w \geq -C$. This proves (\ref{w-lower}) and the proof of Lemma \ref{lemma-c0} is complete. Q.E.D.

\subsection{Second order estimate}
For the second order estimate, we will use the test function
\be
G(x,t) = \log {\rm Tr} \, h - A \varphi + {B\over 2} F^2,
\ee
where $B,A>1$ are constants to be determined. This test function was introduced in \cite{PPZ10} to obtain long-time existence of the Anomaly flow (\ref{af2}) with conformally K\"ahler initial data. It includes the term $F^2$, which is new compared to the standard test function \cite{Yau} used in the proof of the Calabi conjecture. It also differs from the test functions used in previously studied flows with concave operator \cite{Cao,Gill,PT,Gabor,TW1,Zhang}. The key observation is that the differentiation of the $F^2$ term will contribute good quadratic third order terms, which can be used to control the bad terms due to the lack of concavity.

\begin{lemma} \label{lemma-trh}
Let $u$ solve $\p_t u =F(e^{-f} \det \, h)$ with $\chi+i \ddb u>0$ on $X \times [0,T]$. Then
\be
\sup_{X \times [0,T]} {\rm Tr} \, h \leq C,
\ee
where $C$ depends on the initial data, $(X,\chi)$, $f$, $F$.
\end{lemma}

\par {\it Proof:} Combining our previous calculations (\ref{evol-phi}), (\ref{evol-F2}), (\ref{evol-Trh}), we have
\bea \label{evol-G}
(\p_t - L) G 
&=&
 {F' e^{-f}\det \, h \over {\rm Tr} \, h} \bigg\{  {1 \over {\rm Tr} \, h} g^{j \bar{k}} \nabla_{\bar{k}}{\rm Tr} \, h  \nabla_j {\rm Tr} \, h - \chi^{p \bar{q}} g^{j \bar{r}} g^{s \bar{k}} \nabla_p g_{\bar{r} s}  \nabla_{\bar{q}} g_{\bar{k} j}
\nonumber\\
&&
- g^{j \bar{k}} R_{\bar{k} p}{}^r{}_j h^p{}_r +   R  + {1 \over (\det \, h)^2} |\nabla \det \, h|^2 - {2 \over \det \, h}  \Re \langle \nabla f, \nabla \det \, h \rangle \nonumber\\
&&
+  {F'' e^{-f} |\nabla f|^2 \det \, h \over F'} +  {F'' e^{-f} \over F' (\det \, h)} |\nabla \det \, h|^2 \nonumber\\
&&-  2 {F'' e^{-f} \over F'} \Re \langle \nabla f, \nabla \det \, h \rangle +  e^f \Delta e^{-f} \bigg\} 
\nonumber\\
&&
 -AF +A n F' e^{-f} \det \, h -A F' e^{-f} \, \det \, h\, g^{j\bar k} \chi_{\bar{k} j} + {A \over V} \int_X F \, \chi^n
 \nonumber\\
&&
- B (F')^3 e^{-3f} \, \det \, h\, g^{j\bar k} \nabla_j \det \, h  \nabla_{\bar{k}} \det \, h \nonumber\\
&&- B (F')^3 e^{-3f} \, (\det \, h)^3 g^{j\bar k} \nabla_j f \nabla_{\bar{k}} f
 \nonumber\\
&&
+ 2 B (F')^3 e^{-3f} \, (\det \, h)^2 \, \Re \{ g^{j\bar k} \nabla_j f \nabla_{\bar{k}} \det \, h \}.
\eea
Since $\det \, h$ has a uniform positive lower and upper bound, and $F: (0,\infty) \rightarrow {\bf R}$ is smooth and $F'>0$, it follows that
\bea
{1\over C} \leq F'(e^{-f} \det\, h) \leq C, \ \ \ F''(e^{-f} \det\, h) \leq C
\eea
for some constant $C$ depending on $f,F$ and the bounds on $\det \, h$. Then, we can first estimate
\bea
&&{F' e^{-f}(\det \, h) \over {\rm Tr} \, h} 
 \bigg\{ -  g^{j \bar{k}} R_{\bar{k} p}{}^r{}_j h^p{}_r +   R+ {F'' e^{-f} \det\, h \over F' } |\nabla f|^2 + e^f \Delta e^{-f} \bigg\}
\nonumber\\
&\leq &
{C \over {\rm Tr} \, h} \left\{ ({\rm Tr} \, h^{-1} ) \left( {\rm Tr} \, h\right) + 1\right\} \leq C {\rm Tr} \, h^{-1},
\eea
where $C$ is a constant depending on $f, F$ and the curvature of the background metric $\chi$.
\smallskip
\par We also observe that
\bea
{1\over {\rm Tr} \, h}|\nabla \det\, h|^2 = {1\over {\rm Tr} \, h}\chi^{j\bar k} \nabla_j \det\, h \, \nabla_{\bar k} \det\, h \leq g^{j\bar k}\nabla_j \det\, h \, \nabla_{\bar k} \det\, h.
\eea
Using this inequality, we estimate
\bea
&& {F' e^{-f}(\det \, h) \over {\rm Tr} \, h} \bigg\{ {F'' e^{-f}  \over F' (\det \, h)} |\nabla \det \, h |^2 + {1 \over (\det \, h)^2} |\nabla \det \, h|^2 \nonumber\\
&& - \bigg[ {1 \over \det \, h} + {F'' e^{-f} \over F'} \bigg]  2 \Re \langle \nabla f, \nabla \det \, h \rangle\bigg\}  
 \nonumber\\
 &\leq & C \, g^{j\bar k}\nabla_j \det\, h \, \nabla_{\bar k} \det\, h + C \, {\rm Tr} \, h^{-1}
\eea
where $C$ is a constant depending on $f, F$ and the lower and upper bound on $\det \, h$.
\smallskip
\par Next, we estimate 
\bea
-AF +A n F' e^{-f} \det \, h -A F' e^{-f} \, \det \, h\, g^{j\bar k} \chi_{\bar{k} j}+ {A \over V} \int_X F \, \chi^n
\leq - \tau_1 A {\rm Tr} \, h^{-1} + CA,
\eea
where $\tau_1>0$ and $C>0$ depend on $f$, $F$, $n$, $(X,\chi)$, and the bounds on $\det \, h$.
\smallskip
\par For the terms coming from differentiating ${B\over 2} F^2$, we can using the uniform lower and upper bound on $\det \, h$ and obtain the estimate
\bea
&&
- B (F')^3 e^{-3f} \, \det \, h\, g^{j\bar k} \nabla_j \det \, h  \nabla_{\bar{k}} \det \, h 
- B (F')^3 e^{-3f} \, (\det \, h)^3 g^{j\bar k} \nabla_j f \nabla_{\bar{k}} f
 \nonumber\\
&&
+ 2 B (F')^3 e^{-3f} \, (\det \, h)^2 \, \Re \{ g^{j\bar k} \nabla_j f \nabla_{\bar{k}} \det \, h \}
\nonumber\\
&\leq &
- \tau_2 B \, g^{j\bar k} \nabla_j \det \, h  \nabla_{\bar{k}} \det \, h + g^{j\bar k} \nabla_j \det \, h  \nabla_{\bar{k}} \det \, h + C B^2 g^{j\bar k} \nabla_j f \nabla_{\bar{k}} f
\nonumber\\
&\leq &
(1- \tau_2 B) \, g^{j\bar k} \nabla_j \det \, h  \nabla_{\bar{k}} \det \, h + C B^2 \, {\rm Tr} \, h^{-1}.
\eea

\smallskip

Finally, we recall the inequality of Yau and Aubin \cite{Aubin, Yau}, 
\bea \label{aubin-yau}
 {1 \over {\rm Tr} \, h} g^{j \bar{k}} \nabla_{\bar{k}}{\rm Tr} \, h  \nabla_j {\rm Tr} \, h - \chi^{p \bar{q}} g^{j \bar{r}} g^{s \bar{k}} \nabla_p g_{\bar{r} s}  \nabla_{\bar{q}} g_{\bar{k} j} \leq 0.
 \eea
 Indeed, at a point where $\chi_{\bar{k} j} = \delta_{kj}$ and $g_{\bar{k} j} = \lambda_j \delta_{kj}$, then ${\rm Tr} \, h = \sum_p g_{\bar{p} p}$ and
\bea
g^{j \bar{k}} \nabla_{\bar{k}}{\rm Tr} \, h  \nabla_j {\rm Tr} \, h &=& \sum_p {1 \over \lambda_p} \bigg| \sum_i \nabla_p g_{\bar{i} i} \bigg|^2 \nonumber\\
&=& \sum_p {1 \over \lambda_p} \bigg| \sum_i {\nabla_p g_{\bar{i} i} \over \lambda_i^{1/2}} \lambda_i^{1/2} \bigg|^2 \nonumber\\
&\leq& \left( \sum_i \lambda_i \right) \sum_p {1 \over \lambda_p} \sum_i {|\nabla_p g_{\bar{i} i}|^2 \over \lambda_i} \nonumber\\
&=& \left( \sum_i \lambda_i \right) \sum_{i,p} {1 \over \lambda_p\lambda_i }  |\nabla_i g_{\bar{i} p}|^2,
\eea
where we used the Cauchy-Schwarz inequality and the K\"ahler condition $\nabla_i g_{\bar{k} p} = \nabla_p g_{\bar{k} i}$. Thus
\bea
{1 \over {\rm Tr} \, h} g^{j \bar{k}} \nabla_{\bar{k}}{\rm Tr} \, h  \nabla_j {\rm Tr} \, h &\leq& \sum_{i,j,p} {1 \over \lambda_p\lambda_i }  |\nabla_j g_{\bar{i} p}|^2 = \chi^{p \bar{q}} g^{j \bar{r}} g^{s \bar{k}} \nabla_p g_{\bar{r} s}  \nabla_{\bar{q}} g_{\bar{k} j},
\eea
proving (\ref{aubin-yau}). Putting all these estimates into (\ref{evol-G}), we arrive at the following inequality
\bea
&&(\p_t -L) G 
\\
&\leq &
(1+ C - \tau_2 B) \, g^{j\bar k} \nabla_j \det \, h  \nabla_{\bar{k}} \det \, h + \{ C - \tau_1 A + C B^2 \}  \, {\rm Tr} \, h^{-1} + CA. \nonumber
\eea
We can choose $B\gg 1$ such that $1+ C - \tau_2 B \leq 0$. Next, we may choose $A\gg B\gg 1$ such that $C - \tau_1 A + C B^2 \leq -1$. Then, if $G$ attains a maximum on $X \times [0,T]$ at a point $(x_0,t_0)$ with $t_0>0$, by the maximum principle we have the inequality
\be
0 \leq (\p_t - L)G \leq - {\rm Tr} \, h^{-1} + CA.
\ee
It follows that the eigenvalues of $h$ are bounded below at $(x_0,t_0)$. The product of the eigenvalues of $h$ is given by $\det \, h$, which is uniformly bounded along the flow. Thus the eigenvalues of $h$ are bounded above at $(x,t_0)$, and so 
\bea
{\rm Tr} \, h (x_0,t_0) \leq C.
\eea
 Therefore
\be \label{C2-G-est}
G(x,t) \leq G(x_0,t_0) \leq C + A \| \varphi \|_{L^\infty(X \times [0,T])} + {B \over 2} \| F^2 \|_{L^\infty(X \times [0,T])}.
\ee
If $G$ attains a maximum at $t_0=0$, we have already have $G(x,t) \leq C$. Therefore 
\be
\log {\rm Tr} \, h \leq C
\ee
 along the flow. Q.E.D.

\begin{lemma} \label{lemma-ellipticity}
  Let $u$ solve $\p_t u =F(e^{-f} \det \, h)$ with $\chi+i \ddb u>0$ on $X \times [0,T]$. Let $\lambda_i$ denote the eigenvalues of
\be
h^i{}_j = \chi^{i \bar{k}} (\chi_{\bar{k} j} + u_{\bar{k} j}).
\ee
Then for each $i \in \{1, \dots, n \}$, we have the estimate
\be
C^{-1} \leq \lambda_i \leq C, \ \ \|i \ddb u \|_{L^\infty(X \times [0,T])} \leq C,
\ee
where $C$ depends on the initial data, $(X,\chi)$, $f$, $F$.
\end{lemma}

{\it Proof:} The previous lemma gives the upper bound
\be
\lambda_i \leq C.
\ee
By Lemma \ref{lemma-deth}, we know that
\be
C^{-1} \leq \lambda_1 \lambda_2 \cdots \lambda_n \leq C,
\ee
along the flow. If $\lambda_n$ is the smallest eigenvalue, then
\be
\lambda_n = {1 \over \lambda_1 \lambda_2 \cdots \lambda_{n-1}} (\lambda_1 \lambda_2 \cdots \lambda_n) \geq C^{-1},
\ee
establishing the lower bound. Q.E.D.

\

\subsection{Higher order estimates}

In this section, we will derive second order H\"older estimates. We will obtain the following H\"older estimates for $u$
 \bea\label{c2alpha}
\| \p_t u \|_{C^{\delta,\delta/2}(Q)} + \| u_{\bar{k} j} \|_{C^{\delta,\delta/2}(Q)} \leq C,
 \eea
for $0<\delta<1$. Here we use the following notation: on a ball $B_R(0)$ and cylinder $Q= B_R \times (T_0,T)$, for functions $w: B_R \rightarrow {\bf R}$ and $u: Q \rightarrow {\bf R}$, we define
\be
\| w \|_{C^\delta(B_R)} = \| w \|_{L^\infty(B_r)} + \sup_{x \neq y \in B_R} {|w(x)-w(y)| \over |x-y|^\delta},
\ee
and
\be
\| u \|_{C^{\delta,\delta/2}(Q)} = \| u \|_{L^\infty(Q)} + \sup_{(x,t) \neq (y,s) \in Q} {|u(x,t)-u(y,s)| \over (|x-y|+|t-s|^{1/2})^\delta}.
\ee

In general, for real fully nonlinear parabolic equations with concave operator, estimate (\ref{c2alpha}) follows directly from the Krylov $C^{2,\alpha}$ estimate \cite{Krylov82} (or see Theorem 14.7 in \cite{Lieberman}) (which is the parabolic version of the Evans-Krylov theorem for concave elliptic equations). However, the theorem is not applicable here since we do not assume any concavity of our operator $F(e^{-f} \det\, h)$. In \cite{PPZ10}, we followed the argument given by Tsai \cite{Tsai} in the study of the inverse Gauss curvature flow to derive estimate (\ref{c2alpha}) when $F(x)=x$. The proof there also works in the general setting of arbitrary $F$ satisfying $F'>0$. For completeness, we will follow the presentation in \cite{PPZ10} and include details here.

\begin{lemma} \label{c2-alpha-estimate}
Let $u$ be a solution to $\p_t u = F(e^{-f} \det \, h)$ with $\chi+ i \ddb u>0$ on $X \times [0,T]$. Suppose there exists $\Lambda>0$ such that
\be \label{unif-para}
\Lambda^{-1} \chi_{\bar{k} j} \leq g_{\bar{k} j}(x,t) \leq \Lambda \chi_{\bar{k} j},
\ee
and
\be
\bigg\| u - {1 \over V} \int_X u \, \chi^n \bigg\|_{L^\infty(X \times[0,T])} + \| \p_t u \|_{L^\infty(X \times [0,T])} + \| i \ddb u \|_{L^\infty(X \times[0,T])}\leq \Lambda
\ee
along the flow. Let $B_1$ be a coordinate chart on $X$ such that $B_1 \subset {\bf C}^{n}$ is a unit ball. Then there exists $0<\alpha<1$ and $C>0$, depending on $\chi$, $f$, $F$, $T$, $\Lambda$,  such that on $Q=B_{1/2} \times [{T \over 2},T]$,
 \bea
\| \p_t u \|_{C^{\alpha,\alpha/2}(Q)} + \| u_{\bar{k} j} \|_{C^{\alpha,\alpha/2}(Q)} \leq C.
 \eea
\end{lemma}

{\it Proof:} Differentiating the equation in time gives
\be
\p_t F = (F' e^{-f} \det \, h) g^{j \bar{k}} \nabla_j \nabla_{\bar{k}} F.
\ee
By (\ref{unif-para}), this is a linear parabolic PDE for $F$ with uniform ellipticity constant $\Lambda$. By the Krylov-Safonov Harnack inequality \cite{KrSa}, we have the H\"older estimate
\be \label{kry-saf}
\| F(e^{-f} \det \, h) \|_{C^{\alpha,\alpha/2}(Q)} \leq C \| F(e^{-f} \det \, h) \|_{L^\infty(X \times[0,T])} \leq C.
\ee
Let $H = e^{-f} \det \, h$. Then
\be
{|H(x,t)-H(y,s)| \over (|x-y|+|t-s|^{1/2})^\alpha} = {1 \over F'(\theta)} {|F(H(x,t))-F(H(y,s))| \over (|x-y|+|t-s|^{1/2})^\alpha}
\ee
for some $\theta$ between $H(x,t)$ and $H(y,s)$ by the Mean Value Theorem. Since $F'>0$, and $H$ lies in a compact set by (\ref{unif-para}), we obtain
\be
{|H(x,t)-H(y,s)| \over (|x-y|+|t-s|^{1/2})^\alpha} \leq C,
\ee
by (\ref{kry-saf}). Since $e^{-f}$ is a smooth given function, we conclude that
\be
\| \det \, h \|_{C^{\alpha,\alpha/2}(Q)} \leq C.
\ee
In particular, covering $X$ with coordinate balls $B_{1/2}$, we obtain a constant $C$ such that
\be
\| \det \, h (\cdot , t) \|_{C^\alpha(X)} \leq C,
\ee
for all $t \in [T/2,T]$. By the estimate of Y. Wang \cite{Wang} (see also \cite{TWWY}) on complex Monge-Amp\`ere equations, we have 
\be \label{space-c2a}
\| u_{\bar{k} j} (\cdot , t) \|_{C^\delta(X)} \leq C,
\ee
for all $t \in [T/2,T]$, for some $C>0$ and $0<\delta<1$ depending on $(X,\chi)$, $\Lambda$, $\alpha$. 
\bigskip
\par To estimate the full space-time H\"older norm, we next need to consider variations in time of $u_{\bar{k} j}$. For this, we use the argument given in \cite{PPZ10} (see also \cite{Tsai}). Let $\eta>0$. For $t$ and $t+\eta$ in the interval $[0,T]$, and $x \in B_1$, we have
\be
\log \det \, h(x,t) - \log \det \, h(x,t+\eta) = \int_0^1 {d \over ds} \log \det (s h(x,t) + (1-s) h(x,t+\eta)) ds.
\ee
Therefore, since $h^i{}_j = \chi^{i \bar{k}} g_{\bar{k} j}$,
\be \label{translate-time}
\log {\det \, h(x,t) \over \det \, h(x,t+\eta)} = a_\eta^{j \bar{k}} (g_{\bar{k} j}(x,t) - g_{\bar{k} j}(x,t+\eta)),
\ee
where
\be
a_\eta^{j \bar{k}} = \int_0^1 (s h(x,t) + (1-s) h(x,t+\eta))^{-1}{}^j{}_p \chi^{p \bar{k}} \, ds.
\ee
We can write (\ref{translate-time}) as
\be \label{u-eta-pde}
a_\eta^{j \bar{k}} \nabla_j \nabla_{\bar{k}} u_\eta = \psi_\eta,
\ee
where
\be
u_\eta(z,t) = {u(z,t) - u(z,t+\eta) \over |\eta|^{\delta/4}},
\ee
and
\be
\psi_\eta(z,t) = {1 \over |\eta|^{\delta/4}} \log {\det \, h(x,t) \over \det \, h(x,t+\eta)}.
\ee
By the H\"older estimate in space (\ref{space-c2a}) and the uniform ellipticity (\ref{unif-para}), one can show that there exists a $C>0$ independent of $\eta$ such that
\be
C^{-1} \chi^{j \bar{k}} \leq a_\eta^{j \bar{k}}(\cdot, t) \leq C \chi^{j \bar{k}}, \ \ \| a_\eta^{j \bar{k}} (\cdot, t) \|_{C^\delta(B_1)} \leq C,
\ee
and
\be
\| \psi_\eta(\cdot, t) \|_{C^{\delta/4}(B_1)} \leq C,
\ee
for any $t \in [T/2,T-\eta]$. For more details on this step, see Lemma 6 in \cite{PPZ10}. Applying the elliptic Schauder estimates to (\ref{u-eta-pde}) at each time, we obtain
\be
\| u_\eta(\cdot,t) \|_{C^2(B_{1/2})} \leq C(1+\| u_\eta(\cdot,t) \|_{L^\infty(X)}),
\ee
for $t \in [T/2,T-\eta]$. By the estimate $\| \p_t u \|_{L^\infty(X \times [0,T])} \leq \Lambda$, we have
\be
\| u_\eta (\cdot, t) \|_{L^\infty(X)} \leq C,
\ee
for $t \in [T/2,T-\eta]$. Therefore
\be
{|u_{\bar{k} j}(z,t) - u_{\bar{k} j}(z,t+\eta)| \over |\eta|^{\delta/4}} \leq C.
\ee
Combining this with the spacial variation estimate (\ref{space-c2a}), we obtain
\be
\| u_{\bar{k} j} \|_{C^{\delta/2,\delta/4}(Q)} \leq C.
\ee
The estimate on $\p_t u$ follows from the equation $\p_t u = F(e^{-f} \det \, h)$. Q.E.D.

\subsection{Convergence}
Combining our estimates, we have the following result.

\begin{proposition}
Let $(X,\chi)$ be a compact K\"ahler manifold. Let $f \in C^\infty(X,{\bf R})$ and $F: (0,\infty) \rightarrow {\bf R}$ be a smooth function satisfying $F'>0$. Let $u_0$ be a smooth function satisfying $\chi + i \ddb u_0 > 0$. Then there exists a smooth function $u: X \times [0,\infty) \rightarrow {\bf R}$ solving
\be
\p_t u =F(e^{-f} \det \, h), \ \ \chi+i \ddb u>0,
\ee
with $u(x,0)=u_0(x)$. Furthermore, if we let $\varphi = u - {1 \over V} \int_X u \, \chi^n$, then
\be \label{phi-estimates}
\| \varphi \|_{C^k(X\times[0,\infty))} \leq C_k,
\ee
where $C_k$ depends on $u_0$, $(X,\chi)$, $f$, $F$, $k$.
\end{proposition}

{\it Proof:} By standard parabolic theory (see e.g. \cite{Tosatti}), there exists $T>0$ such that a solution exists on $X \times [0,T)$. The differentiated equation is
\be \label{diff-eqn}
{d \over dt} (\nabla_\ell u) = F' e^{-f} (\det \, h) g^{j \bar{k}} \nabla_j \nabla_{\bar{k}} (\nabla_\ell u) - F' e^{-f} (\det \, h) \nabla_\ell f.
\ee
By the estimates obtained so far, this is a linear parabolic PDE for $\nabla_\ell u$ with H\"older continuous, uniformly elliptic coefficients. Interior parabolic Schauder estimates (see e.g. \cite{Kr, Lieberman}) imply that $\nabla_\ell u$ is bounded uniformly on $X \times [T/2,T)$ in the $C^{2+\delta,1+\delta/2}$ norm. The regularity of the coefficients then improves, and a bootstrap argument shows that
\be \label{schauder}
\| \nabla_\ell u \|_{C^k(X \times [T/2,T))} \leq C_k(T).
\ee
Since we have $|\p_t u | \leq C$, then
\be
\| u \|_{C^k(X \times [T/2,T))} \leq C_k(T).
\ee
By the Arzela-Ascoli theorem, if $T<\infty$, then we may extract a subsequence to obtain a smooth function at the final time $u_T$ with $\chi+i \ddb u_T > 0$, and then restart the flow. Thus the flow exists on $[0,\infty)$.
\smallskip
\par The normalized function $\varphi$ satisfies a uniform $C^0$ estimate, independent of time, by Lemma \ref{lemma-c0}. At any point in time $t_0 \geq 1$, the equation (\ref{diff-eqn}) has uniformly elliptic coefficients by Lemma \ref{lemma-ellipticity}, with estimates independent of time. By translating Lemma \ref{c2-alpha-estimate} with $T=1$ to the time interval $[t_0-1/2,t_0+1/2]$ and covering $X$ with coordinate balls $B_{1/2}$, we have uniform H\"older estimates on $u_{\bar{k} j}$ in the time interval $[t_0,t_0+1/2]$. These estimates are independent of $t_0$. Thus the coefficients in the parabolic PDE (\ref{diff-eqn}) have uniform estimates on their H\"older norms. By interior Schauder estimates, all norms of $\nabla_k u = \nabla_k \varphi$ are bounded in the time interval $[t_0,t_0+1/2]$ and are independent of $t_0$. 
\smallskip
\par Since $t_0 \geq 1$ was arbitrary, we obtain a uniform estimate on $\| \varphi \|_{C^k(X \times [1,\infty))}$. The estimate in the small time range $\| \varphi \|_{C^k(X \times [0,1])}$ comes from compactness. Q.E.D.

\bigskip
\par We will prove convergence by using the Harnack inequality, following \cite{Cao} (see also \cite{Gill}). We start by recalling the Harnack inequality.

\begin{proposition}
Let $(M,\hat{g})$ be a compact Hermitian manifold. Let $a(t)$, for $t \in [0,\infty)$, be a family of Hermitian metrics satisfying
  \be
\Lambda \hat{g}^{j \bar{k}} \leq  a^{j \bar{k}}(t) \leq \Lambda \hat{g}^{j \bar{k}}
  \ee
for some $\Lambda >0$. Suppose $v \in C^{\infty}(X \times [0,\infty))$ satisfies $v \geq 0$ and
\be
(\p_t - a^{j \bar{k}} \partial_j \partial_{\bar{k}}) v = 0.
\ee
Then
\be
\sup_{X} v(x,{1/2}) \leq C \inf_{X} v(x,1),
\ee
where $C$ depends only on $(X,\hat{g})$, $a$, $\Lambda$.
\end{proposition}

A full proof of this Harnack inequality in the context of Hermitian geometry, and in the style of Li-Yau \cite{LiYau}, can be found in \cite{Gill}; it also follows from the standard local Harnack inequality for parabolic PDE \cite{Nash,KrSa} and a covering argument similar to the one presented in (\ref{covering1})-(\ref{covering2}).
\bigskip
\par 
Recall that the time-differentiated equation is
\be \label{dot-u-evol}
\p_t \dot{u} = (e^{-f} F' \det \, h) g^{j\bar k} \p_j \p_{\bar{k}} \dot{u},
\ee
where we write $\dot{u} = \p_t u$. Define
\bea
&&\phi_m(x, t) = \sup_{y\in X} \dot{u}(y, m-1) - \dot{u}(x, m-1+t),\nonumber\\
&& \psi_m(x, t) = \dot{u}(x, m-1+t) - \inf_{y\in X} \dot{u}(y, m-1),\nonumber
\eea
where $m$ is a positive integer. These functions satisfy 
\bea
&& \partial_t \phi_m = [e^{-f}F' \det \, h \, g^{j\bar k}](m-1+t) \partial_j \partial_{\bar k} \phi_m,\nonumber\\
&& \partial_t \psi_m = [e^{-f}F' \det \, h \, g^{j\bar k}](m-1+t)  \partial_j \partial_{\bar k} \psi_m.\nonumber
\eea
By the maximum principle, $\phi_m \geq 0$ and $\psi_m \geq 0$. Therefore, we can apply the Harnack inequality to obtain
\bea
\sup_{x\in X} \dot{u}(x, m-1) - \inf_{x\in X} \dot{u}(x, m-1/2) &\leq& C \left( \sup_{x\in X} \dot{u}(x, m-1) - \sup_{x\in X} \dot{u}(x, m)\right)\nonumber\\
\sup_{x\in X} \dot{u}(x, m-1/2) - \inf_{x\in X} \dot{u}(x, m-1) &\leq& C \left( \inf_{x\in X} \dot{u}(x, m) - \inf_{x\in X} \dot{u}(x, m-1)\right). \nonumber
\eea
Let
\be
\omega(t) = \sup_{X} \dot{u}(\cdot, t) - \inf_{X} \dot{u}(\cdot, t).
\ee
Adding the two inequalities gives
\be
\omega(m-1) + \omega(m-1/2) \leq C\left( \omega(m-1) - \omega(m)\right).
\ee
This implies $\omega(m) \leq \delta \omega(m-1)$ with $\delta = {C-1\over C}<1$. Thus
\be
\omega(m) \leq C \delta^m.
\ee
By the evolution of $\dot{u}$ (\ref{dot-u-evol}), we know by the maximum principle that the supremum in space of $\dot{u}$ is decreasing and the infimum of $\dot{u}$ is increasing. Thus $\omega(t)$ is decreasing. For $t \in [m,m+1]$, it follows that
\be
\omega(t) \leq \omega(m) \leq C \delta^{-1} e^{ - (m+1)\log \delta^{-1}} \leq C \delta^{-1} e^{- (\log \delta^{-1}) t}
\ee
and hence
\be
\omega(t) \leq C \, e^{-\eta t}
\ee
with $\eta = \log{1\over \delta}>0$.

Recall the definition of $\varphi = u - {1 \over V} \int_X u \, \chi^n$. It follows that $\int_X \dot{\varphi} \, \chi^n=0$. Therefore, at each time $t$ there exists a point $y\in X$ such that $\dot{\varphi}(y, t)=0$. It follows that for any $(x, t)$,
\bea
\left| \dot{\varphi}(x, t) \right| = \left| \dot{\varphi} (x, t) - \dot{\varphi}(y, t)\right| =\left| \dot{u}(x, t) - \dot{u}(y, t)\right| \leq C \, e^{-\eta t}.
\eea
This implies $\varphi + {C\over \eta} e^{-\eta t}$ is bounded and monotonically decreasing, and hence its limit exists as $t\rightarrow \infty$. Therefore, we can conclude that $\varphi$ converges pointwise to a bounded function $\varphi_{\infty}$. This convergence is in fact smooth by our estimates (\ref{phi-estimates}).
\smallskip
\par Indeed, first we note that $\varphi_\infty$ is smooth, since a subsequence of $\varphi(x,n)$ converges to a smooth function by the Arzela-Ascoli theorem and the estimates (\ref{phi-estimates}). If $\varphi$ does not converge to $\varphi_\infty$ in $C^\infty$, then there exists a sequence $t_m \rightarrow \infty$, an integer $k$, and $\epsilon>0$ such that
\be \label{phi-contradic}
\| \varphi(t_m) - \varphi_\infty \|_{C^k} \geq \epsilon.
\ee
By the Arzela-Ascoli theorem and our estimate (\ref{phi-estimates}), taking a subsequence gives a sequence $\varphi(t_{m_i})$ converging to a smooth function $\psi_\infty$. But $\varphi(t_{m_i})$ converges pointwise to $\varphi_\infty$, hence $\psi_\infty=\varphi_\infty$, contradicting (\ref{phi-contradic}).
\bigskip
\par We now identify the limiting function $\varphi_\infty$. Since $\dot{\varphi} \rightarrow 0$, this function satisfies
\be \label{stationary}
F[\varphi_\infty] = c,
\ee
for a constant $c$, where we used the notation 
\be
F[\varphi_\infty] = F(e^{-f} \det (\delta^i{}_j + \nabla^i \nabla_j \varphi_\infty)).
\ee
Since $F$ is injective, we must have
\be
e^{-f} \det (\delta^i{}_j + \nabla^i \nabla_j \varphi_\infty) = c_0,
\ee
for a constant $c_0$. This constant can be identified by integrating both sides
\be
\int_X (\chi+ i \ddb \varphi_\infty)^n = c_0 \int_X e^f \chi^n.
\ee
By Stokes theorem,
\be
c_0 = {V \over \int_X e^f \chi^n}.
\ee
Therefore, the constant $c$ in (\ref{stationary}) is given by 
\be
c = F \left( \left\{ {1 \over V} \int_X e^f \chi^n \right\}^{-1} \right).
\ee
This completes the proof of the long-time existence and convergence of $\varphi$ along the flow $\p_t u = F(e^{-f} \, \det \, h)$.

\

\bigskip
Department of Mathematics, Harvard University, Cambridge, MA 02138, USA

\smallskip
spicard@math.harvard.edu

\bigskip
Department of Mathematics, University of California, Irvine, CA 92697, USA

\smallskip
xiangwen@math.uci.edu

\end{document}